\input amstex
\documentstyle{amsppt}
\def \al{\alpha}
\def \bt{\beta}
\def \gm{\gamma}

\def \sd #1#2{{#1}^{\scriptstyle{(#2)}}}

\def \QQ{\Bbb {Q}}
\def \NN{\Bbb {N}}
\def \ZZ{\Bbb {Z}}

\def \mm{\Cal M_m^+}
\def \mn{\Cal M_m}

\pagewidth{125mm} 
\pageheight{185mm} 
\parindent=8mm 
\frenchspacing 
\NoBlackBoxes

\topmatter
\title The ring of multisymmetric functions
\endtitle
\author Francesco Vaccarino\endauthor

\address Dipartimento di  Matematica - Universit\'a di Torino - Via Carlo Alberto, 
nr.10 - 10123 - Torino - Italy \endaddress

\curraddr Dipartimento di  Matematica - Universit\'a di Torino - Via Carlo Alberto, 
nr.10 - 10123 - Torino - Italy\endcurraddr 

\email vaccarino\@ syzygie.it \endemail



\abstract 
Let $R$ be a commutative ring and let $n,m$ be two positive integers. Let $A_R(n,m)$ be the polynomial ring in the commuting independent variables $ x_i(j)$ with $i=1, \dots ,m\,;j=1, \dots  ,n$ and coefficients in $R$. 
The symmetric group on $n$ letters $S_n$ acts on $A_R(n,m)$ by means of $\sigma(x_i(j))=x_i(\sigma(j))$ for all $\sigma\in S_n$ and $i=1, \dots ,m\,;j=1, \dots  ,n$. Let us denote by $A_R(n,m)^{S_n}$ the rings of invariants for this action: its elements are usually called multisymmetric functions and they are the usual symmetric functions when $m=1$. In this paper we will give a presentation in terms of generators and relations that holds for any $R$ and any $n,m$ answering in this way to a classical question. I would like to thank M.Brion, C.De Concini and C.Procesi, in alphabetical order, for useful discussions.
\endabstract
\endtopmatter

\document

\head Introduction \endhead
Let $R$ be a commutative ring and let $n,m$ be two positive integers. Let $A_R(n,m)$ be the polynomial ring in the commuting independent variables $ x_i(j)$ with $i=1, \dots ,m\,;j=1, \dots  ,n$ and coefficients in $R$. 
The symmetric group on $n$ letters $S_n$ acts on $A_R(n,m)$ by means of $\sigma(x_i(j))=x_i(\sigma(j))$ for all $\sigma\in S_n$ and $i=1, \dots ,m\,;j=1, \dots  ,n$. Let us denote by $A_R(n,m)^{S_n}$ the rings of invariants for this action: its elements are usually called multisymmetric functions and they are the usual symmetric functions when $m=1$.

\noindent
When $m=1$ then $A_R(n,1)\cong R[x_1,x_2,\dots ,x_n]$ and $R[x_1,x_2,\dots ,x_n]^{S_n}$ is freely generated by the elementary symmetric functions $e_1,\dots,e_n$ that are given by the equality

$$\sum_{k=0}^n t^k e_k:=\prod_{i=1}^n(1+tx_i)\eqno(0.1)
$$
where $e_0=1$ and $t$ is a commuting independent variable (see {\cite{4}}).
Furthermore one has
$$e_k(x_1,\dots,x_n)=\sum_{i_1<i_2<\dots <i_k\leq n}x_{i_1}x_{i_2}\cdots
x_{i_k}\eqno(0.2)$$

\noindent
Let $A_R(m):=R[y_1,\dots,y_m]$, where $y_1,\dots,y_m$ are commuting independent variables, let $f\in A_R(m)$ and define $$f(j):=f(x_1(j),\dots ,x_m(j))\,{\text{ for }}\, 1 \leq j\leq n \eqno(0.3)$$ 
Notice that $f(j)\in A_R(n,m)$ for all $1 \leq j\leq n $. 

\noindent
Define $e_k(f):=e_k(f(1),f(2),\dots,f(n))$ i.e.
$$\sum_{k=0}^n t^k e_k(f):=\prod_{i=1}^n(1+tf(i))\eqno(0.4)
$$

\noindent
Let $\Cal M_m$ be the set of monomials in $A_R(m)$ and $\mm$ the set of those of positive degree.
Let $\mu\in\Cal M_m$  and let $\partial_i(\mu)$ denote the degree of $\mu$ in $y_i$, for all $i=1,\dots ,m$. We set $$\partial(\mu):=(\partial_1(\mu),\dots ,\partial_m(\mu))\eqno(0.5)$$ 
for its multidegree. The total degree of $\mu$ is 
$$l(\mu):=\sum_i\partial_i(\mu).\eqno(0.6)$$ 
These degrees can be extended to $A_R(n,m)$ by $\partial(x_i(j))=\partial(y_i)$ as $S_n$-module. 
If $f$ is homogeneous of total degree $l$ then $e_k(f)$ has total degree $kl$ (for all $k$ and $n$).

\noindent
A monomial $\mu\in \mm$ is called {\it{primitive}} it is not a power of another one.
We denote by $\frak M_m^+$ the set of primitive monomials.

\noindent
We are now able to state the first part of our result.

\proclaim{Generators}
The ring of multisymmetric functions $A_R(n,m)^{S_n}$ is generated by the $e_k(\mu)$ with $\mu\in \frak M_m^+$, $k=1,\dots n$ and $l(e_k(\mu))\leq max(n,n(m-1))$.
If $n=p^s$ is a power of a prime and $R=\ZZ$ or $p\cdot 1_R=0$ then at least one generator has degree equal to $max(n,n(m-1))$. 

\noindent
If $R\supset \QQ$ then $A_R(n,m)^{S_n}$ can also be generated by the $e_1(\mu)$ with $\mu\in \mm$ and $l(\mu)\leq n$. 
\endproclaim

\noindent
Let us now find the relations between these generators that is the second part of our result.

\noindent
Let again $m=1$. The action of $S_n$ on $A_R(n,1)\cong R[x_1,x_2,\dots ,x_n]$ preserves the usual degree. We denote by $\Lambda_{R,n}^k$ the group of
invariants of degree $k$, and 
$$\Lambda_{n,R}:=R[x_1,x_2,\dots ,x_n]^{\scriptstyle{S_n}}=\oplus_{k\geq
0}\Lambda_{n,R}^k .$$ 

\noindent
Let $q_n: R[x_1,x_2,\dots ,x_n]@>>> 
R[x_1,x_2,\dots ,x_{n-1}]$ be given by $x_n\mapsto 0$ and
$x_i\mapsto x_i$, for $i=1,\dots,n-1$. This map sends 
$\Lambda_{n,R}^k$ to $\Lambda_{n-1,R}^k$.

\noindent
Denote by $\Lambda^k_R$ the 
limit of the inverse system obtained in this
way. 

\noindent
The ring $\Lambda_R:=\oplus_{k\geq 0}\Lambda^k_R$ is called the ring of {\it{symmetric functions}} (over $R$).

\noindent
It can be shown {\cite{4}} that $\Lambda_R$ is a free polynomial ring freely generated by the (limit of the) $e_k$, that are given by
$$
\sum_{k=0}^{\infty} t^k e_k:=\prod_{i=1}^{\infty}(1+tx_i). \eqno(0.7)
$$
Furthermore the kernel of the natural projection $\pi_n:\Lambda_R @>>> \Lambda_{n,R}$ is generated by the $e_{n+k}$ with $k\geq n$.
 
\noindent
In a similar way we build a limit of multisymmetric functions. 
For any $a\in \NN^m$ we set $A_R(n,m,a)$ for the linear span of the monomials of multidegree $a$.
One has 
$$A_R(n,m)=\bigoplus_{a\in \NN^m}A_R(n,m,a).\eqno(0.8)$$

\noindent
Let $h$ be an integer with $h>n$. Let $\pi_{n}^h : A_R(h,m)@>>> A_R(n,m)$ be given by
$$\pi_{n}^h(x_i(j))=\cases 0 \;&{\text{if}}\,\, j>n \\
x_i(j) \;&{\text{if}}\,\, j\leq n
\endcases \,\,\,\,\,\;\;\;\;\;\;\;{\text{ for all }}\, i.\eqno(0.9)
$$
then we prove that
$$ \pi_{n}^h(A_R(h,m)^{S_h})=A_R(n,m)^{S_n}\eqno(0.10)$$
\noindent
For any $a \in \NN^m$ set 
$$A_R(\infty,m,a):=\lim_{\leftarrow } A_R(n,m,a)^{S_n} \eqno(0.11)$$ 
where the projective limit is taken with respect to $n$ over the projective system 
$(A_R(n,m,a)^{S_n},\pi_{n,a}^{n+1})$ with 
$\pi_{n,a}^{n+1} : A_R(h,m,a)\rightarrow A_R(n,m,a)$ the restriction of $\pi_{n}^{n+1}$.  

\noindent
Set 
$$A_R(\infty,m):=\bigoplus_{a\in \NN^m} A_R(\infty,m,a).\eqno(0.12)$$

\noindent
We set, by abuse of notation, $$e_j(f):=e_j((f(1),\dots,f(k),\dots))\in A_R(\infty,m) \eqno(0.13)$$ 
with $j\in \NN$ and $f\in A(m)^+$, the augmentation ideal, i.e.
$$
\sum_{k=0}^{\infty} t^k e_k(f):=\prod_{i=1}^{\infty}(1+tf(i)). \eqno(0.14)
$$

\noindent
We can now state the second part of our main result.

\proclaim{Relations}
\roster
\item
The ring $R\otimes A_{\ZZ}(\infty,m)$ is a polynomial ring freely generated by the (limit of) the $1_R\otimes e_k(\mu)$ with $\mu\in\frak M_m^+, \, k\in \NN$. 
The the kernel of the natural projection 
$$\pi_n : R\otimes A_{\ZZ}(\infty,m)@>>> A_R(n,m)^{S_n}$$
is generated by the elements
$$e_{n+k}(f)=0 \,{\text{ with }} \, k\geq 1 \,{\text{ and }}\, f\in A_R(m)^+.$$ 
\item If $R\supset \QQ$ then $A_R(\infty,m)$ can be freely generated by the $e_1(\mu)$ with $\mu\in\mm$. The the kernel of the natural projection is then generated by $e_{n+1}(f)$ with $f\in A_R(m)^+$ .
\endroster
\endproclaim

\head 1. Generators \endhead

\noindent
In order to prove the theorem on generators stated in the introduction we need to introduce some machinery.

\noindent
Let $k\in \NN$, $f_1\dots,f_k\in A_R(m)$ and $t_1,\dots,t_k$ be commuting independent variables, define elements $e_{(\al_1,\dots,\al_k)}(f_1,\dots,f_k)\in A_R(n,m)^{S_n}$ by
$$\sum_{\sum\al_j \leq n}t_1^{\al_1}\cdots t_k^{\al_k}e_{(\al_1,\dots,\al_k)}(f_1,\dots,f_k):=\prod_{i=1}^n(1+t_1 f_1(i) + \dots + t_k f_k(i))\eqno(1.1)$$

\example{Example 1.2}
Let $n=3$ and $f,g\in A_R(m)$ then 
$$e_{(2,1)}(f,g)= f(1)f(2)g(3)+f(1)g(2)f(3)+g(1)f(2)f(3).$$
If $n=4$ then
$$\aligned e_{(2,1)}(f,g) = &f(1)f(2)g(3)+f(1)g(2)f(3)+g(1)f(2)f(3)+\\
														&f(1)f(2)g(4)+f(1)g(2)f(4)+g(1)f(2)f(4)+\\
														&f(1)f(3)g(4)+f(1)g(3)f(4)+g(1)f(3)f(4)+\\
											      &f(2)f(3)g(4)+f(2)g(3)f(4)+g(2)f(3)f(4)\endaligned$$
											      \endexample
											      
\noindent
Let $k=m$ and $f_j=y_j$ for $j=1,\dots,m$ then the $e_{(\al_1,\dots,\al_m)}(y_1,\dots,y_m)$ 
with $\sum\al_j\leq n$ are the well known {\it{elementary multisymmetric function}} that generates $A_R(n,m)^{S_n}$ when $R\supset\QQ$ see {\cite{3}} or {\cite{6}} and are given by
$$\prod_{i=1}^n(1+\sum_{j=1}^m t_jx_j(i))=\sum_{\sum\al_j\leq n}t_1^{\al_1}\cdots t_m^{\al_m}e_{(\al_1,\dots,\al_m)}(y_1,\dots,y_m).\eqno(1.3)$$

\noindent
It can be easily seen that $e_{(\al_1,\dots,\al_k)}(f_1,\dots,f_k)$ is the orbit sum (under the considered action of $S_n$) of $$f_1(1)f_1(2)\cdots f_1(\al_1)f_2(\al_1+1)\cdots f_2(\al_1+\al_2)\cdots f_k(\al_k) \eqno(1.4)$$
It is clear that the $e_{(\al_1,\dots,\al_k)}(f_1,\dots,f_k)=e_{(\al_\tau(1),\dots,\al_\tau(k))}(f_{\tau(1)},\dots,f_{\tau(k)})$ for all $\tau\in S_k$; furthermore if two entries are equal, say $f_1=f_2$, then $$e_{(\al_1,\dots,\al_k)}(f_1,\dots,f_k)=\frac{(\al_1+\al_2)!}{\al_1!\al_2!}e_{(\al_1+\al_2,\dots,\al_k)}(f_1,f_3\dots,f_k).\eqno(1.5)$$
If the arguments $f_i$ of $e_{(\al_1,\dots,\al_k)}(f_1,\dots,f_k)$ are all distinct we say it is {\it{reduced}}.

\noindent
Let $\sd {\NN} {\mn^+}$ be the set of functions $\mn^+ @>>> \NN$ with finite support, we set 
$$\mid \al \mid:=\sum_{\mu\in \mm} \al(\mu)\eqno(1.6)$$
We introduce the elements $e_{\al}(\overline{\mu})\in A_R(n,m)^{S_n}$ with $\mid\al\mid\leq n$ by 
$$\sum_{\mid\al\mid\leq n}t^{\al}e_{\al}(\overline{\mu})=\prod_{i=1}^n(1+\sum_{\mu\in\mn^+}t_{\mu}\mu(i)) \eqno(1.7)$$
where $t_{\mu}$ are commuting independent variables and $t^{\al}:=\prod_{\mu\in\mn^+}t_{\mu}^{\al(\mu)}$ for all
$\al\in \sd {\NN} {\mn^+}$. 
If $\al\in \sd {\NN} {\mn^+}$ is such that $\al(\mu)=k$ and $\al(\nu)=0$ for some $\mu\in\mn^+$ and for all 
$\nu\in\mn^+$ with $\nu\neq\mu$ we see that $e_{\al}(\overline{\mu})=e_k(\mu)$ (see introduction) i.e. the $k$-th elementary symmetric function evaluated at $(\mu(1),\mu(2),\dots,\mu(n))$.

\proclaim{Proposition 1.8}
The (ordered) set 
$$\Cal B_{n,m,R}:=\{e_{\al}(\overline{\mu})\;:\;\mid \al\mid\leq n\}$$ 
is a $R$-basis of $A_R(n,m)^{S_n}$.

\noindent
The (ordered) set $$\Cal B_{n,m,\al,R}:=\{e_{\al}(\overline{\mu})\;:\; \mid\al\mid\leq n {\text{ and }}\partial(e_{\al})(\overline{\mu})=\beta\}$$ 
is a $R$-basis of $A_R(n,m,\beta)^{S_n}$, for all $\beta\in \NN^m$.
\endproclaim
\demo{Proof}
The $e_{\al}(\overline{\mu})$ are complete system of representative (for the action of $S_n$) of the orbit sums of the products
$$\{\mu_1(1)\mu_2(2)\cdots\mu_n(n)\,:\,\mu_i\in\mn \,,\,i=1,\dots,n\}$$ 
furthermore $\partial(e_{\al})(\overline{\mu})=\sum_{\mu\in\mn^+}\al_{\mu}\partial(\mu)$.
\enddemo

\noindent
Let us calculate the product between two elements as above.

\noindent
\proclaim{Proposition 1.9 - Product Formula}

\noindent
Let $k,h\in \NN$, $f_1\dots ,f_k,g_1,\dots , g_h \in A_R(m)$ and $t_1,\dots,t_k, s_1,\dots ,s_h$ be commuting independent variables, set 
$$e_{\al}(\overline{f}):=e_{(\al_1,\dots,\al_k)}(f_1,\dots,f_k) {\text{ and }}e_{\bt}(\overline{g}):=e_{(\bt_1,\dots,\bt_h)}(g_1,\dots,g_h)$$
then
$$e_{\al}(\overline{f})e_{\bt}(\overline{g})=
\sum_{\gm} e_{\gm}(f_1,\dots,f_k,g_1,\dots,g_h,f_1g_1,f_1g_2,\dots,f_kg_1,\dots,f_kg_h)$$
where 
$
\gamma :=(\gamma_{10},\dots,\gamma_{k0},\gamma_{01},\dots,\gamma_{0h},\gamma_{11},\gamma_{12},\dots,\gamma_{k1},\dots,\gamma_{kh})$
is such that
$$
\cases \gamma_{ij}\in \NN \\
\mid \gamma \mid \leq n \\
\sum_{j=0}^h \gamma_{ij}=\al_i \;\;{\text{for}}\;\; i=1,\dots,k\\
\sum_{i=0}^k \gamma_{ij}=\bt_j \;\; {\text{for}}\;\; j=1,\dots,h.
\endcases
$$
\endproclaim
\demo{Proof}
The result follows from
$$\aligned
&\sum_{\sum\al_j,\; \sum\bt_l\leq n}\;\prod_{j=1}^k \;\prod_{l=1}^h t_j^{\al_j}s_l^{\bt_l}e_{\al}(\overline{f})e_{\bt}(\overline{g})= \\
&(\sum_{\sum\al_j\leq n}\prod_{j=1}^kt_j^{\al_j}e_{\al}(\overline{f}))(\sum_{\sum \bt_l\leq n}\prod_{l=1}^h s_l^{\bt_l}e_{\bt}(\overline{g}))= \\
&\prod_{i=1}^n(1+\sum_{j=1}^k t_j f_j(i))\prod_{i=1}^n(1+\sum_{l=1}^h s_l g_l(i))=\\ 
&\prod_{i=1}^n(1+\sum_{j=1}^k t_j f_j(i)+\sum_{l=1}^h s_l g_l(i)+\sum_{j\leq l} t_js_l f_j(i)g_l(i))=\\
&\sum_{\gamma} (\prod_{a=1}^k t_a^{\gamma_{a0}}\prod_{b=1}^h s_b^{\gamma_{0b}}
\prod_{a=1}^k \prod_{b=1}^h (t_a s_b)^{\gamma_{ab}}e_{\gamma}(\overline{f},\overline{g},\overline{fg}))
\endaligned$$
where $\overline{fg}:=(f_1g_1,f_1g_2,\dots,f_kg_1,\dots,f_kg_h)$ and $\gm$ satisfies the condition of the statement.
\enddemo
\example{Example 1.10}
Let us calculate in $A_R(2,3)^{S_2}$
$$e_{(1,1)}(a,b)e_{2}(c)=\sum_{0\leq k,h \leq 1} e_{(1-k,1-h,2-k-h,h,k)}(a,b,c,ac,bc)=e_{(1,1)}(ac,bc)$$ 
since $1-k+1-h+2-k-h+h+k=4-k-h\leq 2$.
\endexample

\noindent
\proclaim{Corollary 1.11}
Let $k\in \NN$, $a_1,\dots ,a_k\in A_R(m)$, $\al=(\al_1,\dots ,\al_k)\in \NN^k$ with $\sum\al_j\leq n$, then $e_{(\al_1,\dots,\al_k)}(a_1,\dots,a_k)$ belongs to the subring of $A_R(n,m)^{S_n}$ generated by the $e_i(\mu)$ with $i=1,\dots,n$ and $\mu$ monomial in the $a_1,\dots,a_k$.
\endproclaim
\demo{Proof}
Let $k, a_1,\dots ,a_k, \al$ be as in the statement, then
$$e_{\al_1}(a_1)e_{(\al_2,\dots,\al_k)}(a_2,\dots,a_k)=$$
$$=e_{(\al_1,\dots,\al_k)}(a_1,\dots,a_k)+\sum e_{\gamma}(a_1,\dots,a_k,a_1a_2,\dots,a_1a_k)$$
where now 
$$\gamma=(\gamma_{10},\gamma_{01},\dots,\gamma_{0h},\gamma_{11},\gamma_{12},\dots,\gamma_{1h})$$
with $h=k-1$ and
$\sum_{j=0}^{h} \gamma_{1j}=\al_1$ with $\sum_{j=1}^h \gamma_{1j}>0$ and $\gamma_{0j}+\gamma_{1j}=\al_j$ for $j=1,\dots,h$.
Thus 
$$\gamma_{10}+\gamma_{01}+\dots+\gamma_{0h}+\gamma_{11}+\dots+\gamma_{1h}=\sum_j\al_j - \sum_{j=1}^h \gamma_{1j}<\sum_j\al_j.$$
Hence  
$$e_{(\al_1,\dots,\al_k)}(a_1,\dots,a_k)=$$ 
$$e_{\al_1}(a_1)e_{(\al_2,\dots,\al_k)}(a_2,\dots,a_k)-\sum e_{\gamma}(a_1,\dots,a_k,a_1a_2,a_1a_3,\dots,a_1a_k)$$
with $\sum_j\al_j < \sum_{r,s}\gm_{rs}.$

\noindent
By induction on $\sum_j\al_j$ it is then possible to write $e_{(\al_1,\dots,\al_k)}(a_1,\dots,a_k)$ 
as a polynomial in $e_j(\mu)$ with $\mu$ monomial in the $a_1,\dots,a_k$.
\enddemo

\example{Example 1.12}
Consider $e_{(2,1)}(a,b)$ in $A_R(3,m)$ as in example 1.2, then 
$$e_{(2,1)}(a,b)=e_2(a)e_1(b)- e_{(1,1)}(a,ab)=e_2(a)e_1(b)-e_1(a)e_1(ab)+e_1(a^2b).$$
\endexample

\noindent
We now recall some basic facts about classical symmetric functions, for further reading on this topic see {\cite{4}}. 

\noindent
In $\Lambda_R$ (see Introduction) we have a further
operation beside the product: the {\it{plethysm}}. Let $g,f\in
\Lambda_R$, we say that $h\in\Lambda_R$ is the plethysm of $g$ by $f$ and
we denote it by $h=g\cdot f$ if $h$ is obtained by substituting the
monomials appearing in $f$ at the place of the variables in $g$.

\noindent
We have another distinguished kind of functions in $\Lambda_R$ beside 
the elementary symmetric ones: the
{\it{powers sums}}. 

\noindent
For any
$r\in \NN$ the $r$-th power sum is 
$$p_r:=\sum_{i\geq 1}x_i^r.$$ 

\noindent
Let $g\in \Lambda_R$,
then the plethysm $g\cdot p_r=g(x_1^r,x_2^r,\dots,x_k^r,\dots)$.
Since the $e_i$ generate $\Lambda_R$ we have that $g\cdot p_r$ can be expressed as a polynomial in the $e_i$, we set  $$P_{h,k}:=e_h\cdot p_k$$
a polynomial in the $e_i$'s.
These result are clearly valid also in $\Lambda_{n,R}$ with the precaution that $e_i=0$ for $i\geq n+1$.

\proclaim{Proposition 1.13}
For all $a\in A_R(m)$, and $k,i\in \NN$, $e_i(a^k)$ belongs
to the subring of $A_R(n,m)^{S_n}$ generated by the $e_j(a)$.
\endproclaim
\demo{Proof} 
Let $a\in A(m)$ and consider $e_h(a^k)\in A(n,m)^{S_n}$, we have (see Introduction)
$$e_h(a^k)=e_h(a(1)^k,\dots,a(n)^k)=P_{h,k}(e_1(a(1),\dots,a(n)),\dots,e_n(a(1),\dots,a(n)))$$
and the result is proved.\enddemo

\noindent
Recall that a monomial $\mu\in \mn^+$ is called {\it{primitive}} if it is not a power of another one and we denote by $\frak M_m^+$ the set of primitive monomials. Let we state again the first part of the main result
\proclaim{Theorem 1.14 - Generators}

\noindent
The ring of multisymmetric functions $A_R(n,m)^{S_n}$ is generated by the $e_k(\mu)$ with $\mu\in \frak M_m^+$, $k=1,\dots n$ and $l(e_k(\mu))\leq max(n,n(m-1))$.
If $n=p^s$ is a power of a prime and $R=\ZZ$ or $p\cdot 1_R=0$ then at least one generator has degree equal to $max(n,n(m-1))$. 

\noindent
If $R\supset \QQ$ then $A_R(n,m)^{S_n}$ can also be generated by the $e_1(\mu)$ with $\mu\in \mm$ and $l(\mu)\leq n$. 
\endproclaim\demo{Proof}
The elements $e_{\al}(\overline{\mu})\in\Cal B_{n,m,R}$, that form a $R$-basis by 1.8, can be expressed as polynomial in $e_i(\mu)$ with $i=1,\dots, n$ and $\mu \in \mm$, by 1.11. If $\mu=\nu^k$, with $\nu\in \frak M_m^+$ then $e_i(\mu)$ can be expressed as a polynomial in the $e_j(\nu)$ in $A(n,m)^{S_n}$, by 1.13. Since for all $\mu\in\mm$ there exist $k\in\NN$ and $\nu\in\frak M_m^+$ such that $\mu=\nu^k$ we have that $A(n,m)^{S_n}$ is generated as commutative ring by the $e_j(\nu)$, with $\nu\in \frak M_m^+$ and $j=1,\dots ,n$.

\noindent
The result then follows by the following due to Fleischmann {\cite{2}}: $A_R(n,m)^{S_n}$ can be generated by elements of total degree $\ell\leq max(n,n(m-1))$, for any commutative ring $R$, with sharp bound if $n=p^s$ a power of a prime and $R=\ZZ$ or $p\cdot 1_R=0$.
If $R\supset \QQ$ then the result follows from Newton's Formulas and the well known result of H.Weyl (see {\cite{3}},{\cite{6}}).
\enddemo

\head 2. Relations \endhead
\noindent
Let now $h$ be an integer with $h>n$. Let $\pi_{n}^h : A_R(h,m)@>>> A_R(n,m)$ be given by
$$\pi_{n}^h(x_i(j))=\cases 0 \;&{\text{if}}\,\, j>n \\
x_i(j) \;&{\text{if}}\,\, j\leq n
\endcases \,\,\,\,\,\;\;\;\;\;\;\;{\text{ for all }}\, i.\eqno(2.1)
$$
this is a multigraded ring epimorphism such that 
$$\prod_{i=1}^h(1+t_1a_1(i)+\cdots +t_ka_k(i))\mapsto \prod_{i=1}^{n}(1+t_1a_1(i)+\cdots +t_ka_k(i))\eqno(2.2)$$ 
where $k, a_j,t_j$ are as usual. 
Hence 
$$\pi_{n}^h(e_{\al}(\overline{\mu}))=\cases e_{\al}(\overline{\mu}) &\;{\text{if}}\; \mid \al \mid \leq n\\           
																						0 &\;{\text{if}}\;\mid \al \mid >n. 
\endcases
\eqno(2.3)$$ 
thus, by Prop.1.8, for all $a\in\NN^m$
the restriction
$$\pi_{n,a}^{h} : A_R(h,m,a)\rightarrow A_R(n,m,a)\eqno(2.4)$$
is such that 
$$ \pi_{h,a}^h(A_R(h,m,a)^{S_h})=A_R(n,m,a)^{S_n}.\eqno(2.5)$$
and then $(A_R(n,m,a)^{S_n},\pi_{n,a}^{n+1})$ is a projective sytem.

\noindent
For any $a \in \NN^m$ set 
$$A_R(\infty,m,a):=\lim_{\leftarrow } A_R(n,m,a)^{S_n} \eqno(2.6)$$ 
where the projective limit is taken with respect to $n$ over the above projective system and
set 
$$\pi_{n,a}:A_R(\infty,m,a)@>>> A_R(n,m,a)^{S_n}\eqno(2.7)$$ 
for the natural projection.

\noindent
Set 
$$A_R(\infty,m):=\bigoplus_{a\in \NN^m} A_R(\infty,m,a).\eqno(2.8)$$
and
$$\pi_n:=\oplus_{a\in\NN^m}\pi_{n,a}.\eqno(2.9)$$

\noindent
Similarly to the classic case ($m=1$) we make an abuse of notation and set $e_{\al}(\overline{\mu}):=\underset{\leftarrow }\to{\lim}\; e_{\al}(\overline{\mu})$, for any $\al \in \sd {\NN}{\mm}$ where $a=\partial(e_{\al}(\overline{\mu}))$. In the same way we set $e_j(f):=\underset{\leftarrow }\to{\lim}\;e_j(f)$ with $j\in \NN$ where $f\in A_R(m)^+$ is homogeneous of positive multidegree so that $j\;\partial(f)=a$. 

\proclaim{Proposition 2.10}Let $a\in \NN^m$. 
\roster
\item The $R$-module $\ker \pi_{n,a}$ is the linear span of 
$$\{e_{\al}(\overline{\mu})\in A_R(\infty,m,a)\,:\,\mid \al \mid >n\}.$$
\item The $R$-module homomorphisms $\pi_{n,a}:A_R(\infty,m,a)@>>> A_R(n,m,a)^{S_n}$ are onto for all $n\in\NN$ and $A_R(\infty,m,a) \cong A_R(n,m,a)^{S_n}$ for all $n \geq \mid a \mid.$
\item The $R$-module $A_R(\infty,m,a)$ is a free $R$-module with basis $$\{e_{\al}(\overline{\mu}) \, : \, \partial(e_{\al}(\overline{\mu}))=a \},$$
\item The $R$-module $A_R(\infty,m)$ is a free $R$-module with basis $$\{e_{\al}(\overline{\mu}) \, : \, \al\in\sd \NN {\mm} \},$$
\endroster
\endproclaim
\demo{Proof}
\roster
\item By 2.5, for all $h>n$ and $a\in\NN^m$, the following is an exact sequence of $R$-modules 
$$0 @>>> \ker \pi_{n,a}^h @>>> A(h,m,a)^{S_h} @>\pi_{n,a}^h>> A(n,m,a)^{S_n} @>>> 0.$$
that splits by (2.3) and Prop.1.8 and the claim follows.
\item If $\sum_{j=1}^m a_j \leq k $ then $\ker \pi_{n,a}^k= 0$ for all $k>n$, indeed $$\partial(e_{\al}(\overline{\mu})=\sum_{\mu \in \mm} \al_{\mu}\;\partial(\mu)=a \Longrightarrow \mid \al \mid \leq \sum_{j=1}^m a_j  \leq n.$$
Hence $A(h,m,a)^{S_h}\cong A(\sum_{j=1}^m a_j,m,a)^{S_h}$ for all $h\geq \sum_{j=1}^m a_j$ and the thesis follows by (2.5).
\item It follows from (1) and (2).
\item It follows form (3) and (2.8)
\endroster
\enddemo

\noindent
\remark{Remark 2.11}
Notice that $A_R(m)^{\otimes n}\cong A_R(n,m)$ as multigraded $S_n$-algebras by means of  
$$f_1\otimes\cdots \otimes f_n \leftrightarrow f_1(1)f_2(2)\cdots f_n(n)\eqno(2.12)$$
for all $f_1,\dots ,f_n \in A_R(m)$.
Hence $A_R(n,m)^{S_n}\cong TS^n(A_R(m))$, where $TS^n(\;-\;)$ denotes the symmetric tensors functor. 
Since $TS^n(A_R(m))\cong R\bigotimes TS^n(A_{\ZZ}(m))$ (see {\cite{1}})we have 
$$A_R(n,m)^{S_n}\cong R\otimes A_{\ZZ}(n,m)^{S_n}\eqno(2.13)$$ 
for any commutative ring $R$.
\endremark

\noindent
We then work with $R=\ZZ$ and we suppress the ${\ZZ}$ subscript for sake of simplicity. 

\remark {Remark 2.14} The $\ZZ$-module $A(\infty,m)$ can be endowed with a structure of $\NN^m$-graded ring such that the $\pi_n$ are $\NN^m$-graded ring homomorphisms: the product $e_{\al}(\overline{\mu})e_{\bt}(\overline{\nu})$, where $\al ,\bt \in \sd {\NN} {\mm}$, is defined using the product formula of Prop.1.9 with no limits on the maximum value of $\mid\gamma\mid$ with $\gamma$ appearing in the summation.
\endremark

\proclaim{Proposition 2.15}
Consider the free polynomial ring 
$$C(m):=\bigoplus_{a\in \NN^m}C(m,a):=\ZZ[e_{i,\mu}]_{i\in \NN,\mu \in \frak M_m^+}$$ 
with multidegree given by $\partial(e_{i,\mu})=\partial(\mu)i$. 

\noindent
The multigraded ring homomorphism
$$\sigma_m:\ZZ[e_{i,\mu}]_{i\in \NN,\mu \in \frak M_m^+}@>>>A(\infty,m)$$
given by
$$\sigma_m: e_{i,\mu} \mapsto e_i(\mu), \,\,{\text{for all}}\,\,i\in \NN,\mu \in \frak M_m^+$$
is an isomorphism, i.e. $A(\infty,m)$ is freely generated as commutative ring by the $e_i(\mu)$ with $i\in \NN,\mu \in \frak M_m^+.$
\endproclaim
\demo{Proof}
Since we defined the product in $A(\infty,m)$ as in Prop.1.9 it is easy to verify, repeating the same reasoning of the previous chapter, that $A(\infty,m)$ is generated as commutative ring by the $e_i(\mu)$ with $i\in \NN,\mu \in \frak M_m^+.$
Hence $\sigma_m$ is onto for all $m\in\NN$.

\noindent
Let $a\in \NN^m$ and consider the restriction $\sigma_{m,a}:C(m,a)@>>>A(\infty,m,a)$.
It is onto as we have just seen. A $\ZZ$-basis of $C(m,a)$ is 
$$\{\prod_{i\in\NN,k\in\NN,\mu \in \frak M_m^+}e_{i,\mu}\,:\,
\sum_{i\in\NN,k\in\NN,\mu \in \frak M_m^+}i\;k\;\partial(\mu)\;=\;a\}$$
on the other hand a $\ZZ$-basis of $A(\infty,m,a)$ is 
$$\{e_{\al}\,:\,
\sum_{\al_{\mu}\in \NN,\mu \in \mm}\al_{\mu}\;\partial(\mu)\;=\;a\}.$$

\noindent
Let $\mu\in\mm$ then there is an unique $k\in\NN$ and an unique $\nu\in\frak M_m^+$ such that $\mu=\nu^k$, hence 
$$\sum_{\al_{\mu}\in \NN,\mu \in \mm}\al_{\mu}\;\partial(\mu)= \sum_{k\in\NN, \al_{\mu}\in \NN, \nu \in \frak M_m^+}\al_{\mu}\;k\;\partial(\nu)$$
so that $C(m,a)$ and $A(\infty,m,a)$ have the same (finite) $\ZZ$-rank and thus are isomorphic. \enddemo

\proclaim{Corollary 1.16}
Let $R\supset\QQ$ then $R\otimes A(\infty,m)$ is a polynomial ring freely generated by the $e_1(\mu)$ with $\mu\in\mm$.
\endproclaim
\demo{Proof}
By Prop.2.15 and Th.1.14.
\enddemo

\proclaim{Theorem 2.17 - Relations}
\roster
\item Let $R$ be a commutative ring: the following is a presentation of $A_R(n,m)^{S_n}$ in terms of generators and relations in the category of $\NN^m$-graded $R$-algebras
$$0 @>>> < \{1_R \otimes e_k(f)\,:\, k>n \;, \; f\in A_R(m)^+\}> @>>> R\otimes A(\infty,m) @>\pi_n>> A_R(n,m)^{S_n}@>>> 0.$$
\item
If $R\supset \QQ$ then $A_R(\infty,m)$ can be freely generated by the $e_1(\mu)$ with $\mu\in\mm$. The the kernel of the natural projection is then generated by $e_{n+1}(f)$ with $f\in A_R(m)^+$ .
\endroster
\endproclaim
\demo{Proof} 
\roster
\item As before we set $R=\ZZ$ and the result will follow by Remark 2.11 and Prop.2.15. 

\noindent
By Prop.2.10 the kernel of $A(\infty,m) @>\pi_n>> A(n,m)^{S_n}$ has basis $\{e_{\al}(\overline{\mu})\,\:\,\mid \al \mid >n \}$, let $V_k$ be the submodule of $A(\infty,m)$ with basis $\{e_{\al}(\overline{\mu})\,\:\,\mid \al \mid=k \}$. Let $A_k$ be the sub-$\ZZ$-module of $\QQ\otimes V_k$ generated by the $e_k(f)$ with $f\in A(m)^+$. Let $g:\QQ\otimes V_k @>>> \QQ$ be a linear form identically zero on $A_k$, then
$$0=g(e_k(f))=g(e_k(\sum_{\mu\in \mm}\lambda_\mu \;\mu))=(\sum_{\mid \al \mid=k}(\prod_{\mu\in \mm}\lambda^{{\al}_{\mu}}_{\mu})\;g(e_{\al}(\overline{\mu}))),$$
for all $\sum_{\mu\in \mm}\lambda_{\mu}\;\mu \in A(m)^+$, hence $g(e_{\al}(\overline{\mu}))=0$ for all $e_{\al}(\overline{\mu})$ with $\mid \al \mid=k$ , thus $g=0$ and $A_k=V_k$. 
\item If $R\supset\QQ$ the result follows from Newton's formulas (1) and Cor.2.16
\endroster
\enddemo

\refstyle{C}


\Refs

\ref \key 1 \by N.Bourbaky \book Elements of mathematics - Algebra II Chapters 4-7 \publ Springer-Verlag \publaddr Berlin \yr 1988 \endref

\ref \key 2 \by P.Fleischmann \paper A new degree bound for vector invariants of symmetric groups \jour Trans. Am. Math. Soc. \vol 350 \pages 1703-1712  \yr 1998 \endref

\ref \no 3\by I.Gelfand,M.Kapranov, A.Zelevinsky \book Discriminants, resultants and  multidimensional determinants \publ Birkahuser  \publaddr Boston \yr 1994 \endref

\ref \no 4\by I.G.Macdonald \book Symmetric Functions and Hall Polynomials - second edition \publ Oxford mathematical monograph \yr 1995 \endref

\ref \no 6\by H.Weyl \book The classical groups\publ Princeton University Press \publaddr Princeton N.J.\yr 1946\endref

\endRefs
\enddocument